\documentclass[11pt]{article}

\usepackage[margin=1in]{geometry}

\usepackage{amsmath,amssymb,amsfonts}
\usepackage{amsthm}

\newtheorem{theo}{Theorem}
\newtheorem{lem}{Lemma}
\newtheorem{cor}{Corollary}

\theoremstyle{definition}
\newtheorem{defn}{Definition}
\newtheorem{ex}{Example}

\theoremstyle{remark}

\newcommand{\N}{\mathbb{N}}
\newcommand{\C}{\mathbb{C}}
\newcommand{\T}{\mathbb{T}}
\newcommand{\Dbar}{\overline{\mathbb{D}}}
\newcommand{\norm}[1]{\lVert#1\rVert}
\newcommand{\abs}[1]{\lvert#1\rvert}
\newcommand{\sm}{\mathcal{S}}

\title{Attractivity Implies Lyapunov Stability for\\
Multidimensional Fornasini--Marchesini Systems}

\author{A. Bourget and T. McMillen\\
Department of Mathematics\\
California State University, Fullerton}

\begin{document}
\maketitle

\begin{abstract}
We study the linear multidimensional Fornasini--Marchesini recurrence with
uniformly bounded boundary data prescribed on the open coordinate faces.
Attractivity is understood in the sense of Pringsheim convergence, in which
all discrete indices tend to infinity. We establish that attractivity alone
implies Lyapunov stability, and consequently is equivalent to asymptotic
stability, for an arbitrary number of discrete directions. Our proof combines
inheritance of attractivity by lower-dimensional subsystems, a Baire-category
uniform boundedness argument, and word-sum estimates to control coordinate
slabs outside the Pringsheim tail. We use a generalized spectral stability
criterion proved in a separate companion preprint, without reproving that
criterion here. We also obtain a uniform exponential interior estimate and
contrast our boundary/convergence convention with an earlier two-dimensional
attractivity theorem.
\end{abstract}

\section{Introduction}
The familiar implication from attractivity to Lyapunov stability for
finite-dimensional linear difference equations is less immediate for
multidimensional systems. For a system indexed by several independent
integers, the ``initial condition'' is an array of boundary states. Its
stability properties therefore depend on the admissible boundary space and
on the interpretation of convergence at infinity.

Fornasini--Marchesini (FM) systems form a classical framework for
multidimensional discrete dynamics and filtering \cite{FM1978,FM1980,Bose}.
A stability analysis for doubly indexed FM systems with bounded boundary
data was developed in \cite{BourgetMcMillenSIMAX}.
Bachelier, David, and the Yeganefars proved that attractivity implies
asymptotic stability for two-dimensional FM and Roesser systems
\cite{BachelierECC}. In their framework the boundary values tend to zero,
and solutions converge as the \emph{sum} of the indices tends to infinity.
Here we allow \emph{all bounded} boundary data, but require convergence only
when \emph{every} index diverges (Pringsheim convergence). These are
mathematically different conventions. A simple scalar example in
Section~\ref{sec:comparison} makes the distinction explicit.

Our main result is that, in this bounded-boundary/Pringsheim framework,
\begin{equation}\label{eq:aim}
\text{attractivity}\quad\Longrightarrow\quad
\text{Lyapunov stability}
\end{equation}
for $m$-dimensional FM models with arbitrary $m\ge 1$ and without
commutativity assumptions on the coefficient matrices. The main obstacle
is that Pringsheim convergence yields no direct bound along coordinate
slabs where one index remains small and the others become arbitrarily large.

A companion preprint \cite{BourgetMcMillenPreprint}  establishes the generalized FM spectral criterion used below. The principal contribution of the present paper is implication \eqref{eq:aim}, which we prove by combining dimensional induction, a Baire-category argument, and uniform estimates on coordinate slabs.

\section{Model and a spectral preliminary}\label{sec:prelim}

Let $A_1,\ldots,A_m\in M_n(\C)$ and
$\alpha=(\alpha_1,\ldots,\alpha_m)\in\N^m$,
where $\N=\{1,2,\ldots\}$ and $\N_0=\N\cup\{0\}$.
The $m$-dimensional Fornasini--Marchesini (FM) system
$\sm(A_1,\ldots,A_m)$ is defined by
\begin{equation}\label{eq:FM}
	x_\alpha=\sum_{k=1}^m A_kx_{\alpha-e_k},
	\qquad \alpha\in\N^m,
\end{equation}
where $x_\alpha\in\C^n$ and $e_k$ denotes the
$k$th standard basis vector of $\mathbb C^m$.
Boundary data are prescribed on
\begin{equation}\label{eq:faces}
	\mathcal H_m=\{\beta\in\N_0^m:\text{exactly one }\beta_k=0\},
\end{equation}
the union of the open coordinate faces.
Let $X_m=\ell^\infty(\mathcal H_m;\C^n)$, with norm
$\norm{\xi}_\infty=\sup_{\beta\in \mathcal H_m}\norm{\xi_\beta}$.
Each $\xi\in X_m$ uniquely determines a solution
of \eqref{eq:FM}.

\begin{defn}
	The system $\sm(A_1,\ldots,A_m)$ is
	\emph{Lyapunov stable} if there exists $C>0$
	such that every $\xi\in X_m$ yields
	$\sup_{\alpha\in\N^m}\norm{x_\alpha}
	\le C\norm{\xi}_\infty$.
	It is \emph{attractive} if, for every $\xi\in X_m$,
	$x_\alpha\to0$ as
	$\min_{1\le k\le m}\alpha_k\to\infty$
	(Pringsheim convergence).
	It is \emph{asymptotically stable} if it is
	both Lyapunov stable and attractive.
\end{defn}

To express solutions, let $W_\beta$ be the sum of all ordered words in
$A_1,\ldots,A_m$ containing exactly $\beta_k$ copies of $A_k$, with
$W_0=I$ and $W_\beta=0$ if some component of $\beta$ is negative.
For $\alpha\in\N^m$, set
\begin{align}\label{eq:R}
R_k(\alpha)=\{\beta\in\N_0^m:\;&\beta_k=\alpha_k-1,\\[-2pt]
 &0\le\beta_j\le\alpha_j-1\ (j\ne k)\}.\notag
\end{align}
Unrolling \eqref{eq:FM} until a boundary face is first reached gives
\begin{equation}\label{eq:representation}
 x_\alpha=\sum_{k=1}^m\sum_{\beta\in R_k(\alpha)}
 W_\beta A_k\xi_{\alpha-\beta-e_k}.
\end{equation}
This word-sum representation, whose general form is proved in
\cite{BourgetMcMillenPreprint}, also follows by induction on
$\abs{\alpha}-m$, where $\abs\alpha=\sum_k\alpha_k$.

Let $\mathbb{T}=\{z\in \mathbb{C} : |z|=1\}$ denote the unit circle. 
The only spectral result needed below is the following theorem from
\cite{BourgetMcMillenPreprint}. 
\begin{theo}[Generalized FM criterion
\cite{BourgetMcMillenPreprint}]\label{thm:FMcriterion}
For $A_k\in M_n(\C)$, the system $\sm(A_1,\ldots,A_m)$ is
asymptotically stable if and only if
\begin{equation}\label{eq:toral}
\rho\!\left(\sum_{k=1}^m z_kA_k\right)<1
\quad\text{for all }(z_1,\ldots,z_m)\in\T^m.
\end{equation}
Equivalently,
$\det(I-\sum_k z_kA_k)\ne0$ on the closed polydisk $\Dbar^m$.
\end{theo}

We isolate the consequence used in the subsequent induction.
\begin{lem}[Summability of word sums]\label{lem:summability}
If $\sm(A_1,\ldots,A_m)$ is asymptotically stable, then
\begin{equation}\label{eq:abs-sum}
 \sum_{\beta\in\N_0^m}\norm{W_\beta}<\infty.
\end{equation}
More precisely, some $r>1$ and $M>0$ satisfy
$\norm{W_\beta}\le Mr^{-\abs\beta}$ for all $\beta\in\N_0^m$.
\end{lem}
\begin{proof}
By Theorem~\ref{thm:FMcriterion}, the function
$F(z)=(I-\sum_kz_kA_k)^{-1}$ is holomorphic on a neighborhood of
$\Dbar^m$. Hence it is holomorphic on a polydisk $\abs{z_k}<R$
for some $R>1$. For $z$ near the origin, the Neumann series gives
\begin{equation}\label{eq:gen}
 F(z)=\sum_{p\ge0}\left(\sum_kz_kA_k\right)^p
     =\sum_{\beta\in\N_0^m}W_\beta z^\beta.
\end{equation}
By uniqueness, the $W_\beta$ are the Taylor coefficients of $F$.
Multivariable Cauchy estimates on the torus $\abs{z_k}=r$, with
$1<r<R$, imply $\norm{W_\beta}\le Mr^{-\abs\beta}$.
Summing the geometric bounds proves \eqref{eq:abs-sum}.
\end{proof}

\section{Attractivity and lower-dimensional systems}
We denote by 
$\sm_i=\sm(A_1,\ldots,A_{i-1},A_{i+1},\ldots,A_m)$ the
$(m-1)$-dimensional system obtained by removing $A_i$.

\begin{lem}[Inheritance of attractivity]\label{lem:inherit}
	If $m\ge2$ and $\sm(A_1,\ldots,A_m)$ is attractive,
	then each $\sm_i$ is attractive.
\end{lem}

\begin{proof}
	By symmetry take $i=m$ and set $\widehat{\sm}=\sm_m$.
	If $\widehat{\sm}$ is not attractive, some bounded datum
	$\widehat\xi$ produces a solution $\widehat x$ with
	$\norm{\widehat x_{\gamma^{(s)}}}\ge\varepsilon>0$
	at points $\gamma^{(s)}\in\N^{m-1}$ tending to
	infinity componentwise.
	
	We construct a single bounded boundary datum for the
	full system by placing copies of $\widehat\xi$ on
	successively distant layers. For $q\ge1$, define
	the boundary layer
	\[
	L_q=\{(\beta,q):\beta\in \mathcal H_{m-1}\}\subset \mathcal H_m.
	\]
	Choose integers $1\le q_1<q_2<\cdots$ tending to
	infinity and prescribe
	$\xi_{(\beta,q_s)}=\lambda_s\widehat\xi_\beta$
	for $\beta\in \mathcal H_{m-1}$, where
	$\lambda_s\in\{-1,1\}$, and set $\xi=0$
	on the remaining boundary points.
	Clearly, $\norm{\xi}_\infty=\norm{\widehat\xi}_\infty$.
	
	Consider $\alpha^{(s)}=(\gamma^{(s)},q_s)$.
	By \eqref{eq:representation}, only boundary layers
	with final coordinate at most $q_s$ can contribute.
	Moreover, contributions from $L_{q_s}$ require
	$k\ne m$ and $\beta_m=0$, so the corresponding
	word sums contain no $A_m$ factors and coincide
	with those of the reduced system.
	 Such paths are precisely those
	of the reduced system, so
	\[
	x_{\alpha^{(s)}}=
	v_s+\lambda_s\widehat x_{\gamma^{(s)}},
	\]
	where $v_s$ depends only on the previously chosen
	signs $\lambda_1,\ldots,\lambda_{s-1}$.
	Since
	$2\norm{w}\le\norm{v+w}+\norm{v-w}$,
	we may choose $\lambda_s\in\{-1,1\}$ recursively
	so that
	\[
	\norm{x_{\alpha^{(s)}}}
	\ge\norm{\widehat x_{\gamma^{(s)}}}
	\ge\varepsilon.
	\]
	As $\alpha^{(s)}$ tends to infinity componentwise,
	this contradicts attractivity.
\end{proof}

\begin{lem}[Uniform tail bound]\label{lem:tail}
If $\sm(A_1,\ldots,A_m)$ is attractive, then there exist
$N\in\N$ and $C_0>0$ such that
\begin{equation}\label{eq:tail}
 \norm{x_\alpha}\le C_0\norm{\xi}_\infty
 \quad\text{whenever }\alpha_k\ge N\ (1\le k\le m).
\end{equation}
\end{lem}

\begin{proof}
Let $T_\alpha:X_m\to\C^n$ assign to $\xi$ the solution value
$x_\alpha$. Each $T_\alpha$ is a bounded linear operator by
\eqref{eq:representation}. Linearity and attractivity give
$T_\alpha\xi\to0$ in the Pringsheim sense for every $\xi\in X_m$.
For integers $N,j\ge1$, put
\[
 E_{N,j}=\left\{\xi\in X_m:
 \sup_{\alpha_1,\ldots,\alpha_m\ge N}
 \norm{T_\alpha\xi}\le j\right\}.
\]
Each $E_{N,j}$ is closed and $X_m=\bigcup_{N,j\ge1}E_{N,j}$.
Since $X_m$ is a Banach space, Baire's theorem gives a ball
$B(\xi_0,\delta)$ contained in some $E_{N,j}$.
For $\norm\eta_\infty<\delta$, both $\xi_0$ and $\xi_0+\eta$
belong to this set. Therefore
$\norm{T_\alpha\eta}\le2j$ on the tail.
Homogeneity implies
$\sup_{\alpha_k\ge N}\norm{T_\alpha}\le2j/\delta$,
which proves \eqref{eq:tail}.
\end{proof}

\section{Main result}\label{sec:main}
\begin{theo}\label{thm:main}
For the multidimensional FM system \eqref{eq:FM}, with arbitrary bounded
boundary data on $\mathcal H_m$, attractivity in the Pringsheim sense implies
Lyapunov stability for every $m\ge1$.
\end{theo}

\begin{proof}
We proceed by induction on $m$. For $m=1$, attractivity of
$x_k=A_1x_{k-1}$ implies $\rho(A_1)<1$ and hence Lyapunov stability.
Assume the assertion for $(m-1)$-dimensional systems, and suppose
$\sm(A_1,\ldots,A_m)$ is attractive. Lemma~\ref{lem:tail}
provides constants $N$ and $C_0$ such that
\eqref{eq:tail} holds when all $\alpha_k\ge N$.
We need only bound the finitely many coordinate slabs
$\{\alpha\in\N^m:\alpha_i<N\}$.

Fix $i$. By Lemma~\ref{lem:inherit}, the lower-dimensional system
$\sm_i$ is attractive. By the induction hypothesis it is Lyapunov
stable, hence asymptotically stable. Applying
Lemma~\ref{lem:summability} to the word sums
$W^{(i)}_\gamma$ of the matrices excluding $A_i$, we obtain
\begin{equation}\label{eq:Li}
 L_i:=\sum_{\gamma\in\N_0^{m-1}}
            \norm{W^{(i)}_\gamma}<\infty.
\end{equation}

A word containing exactly $p$ copies of $A_i$ has a unique
factorization $U_0A_iU_1A_i\cdots A_iU_p$, where each $U_j$ is a
(possibly empty) word in the other matrices. Grouping terms by the
multiplicities $\gamma^j$ of those words gives, with the $i$th
coordinate separated from the others,
\begin{equation}\label{eq:word-decomp}
 W_{(p,\gamma)}=
 \sum_{\gamma^0+\cdots+\gamma^p=\gamma}
 W^{(i)}_{\gamma^0}A_iW^{(i)}_{\gamma^1}
 \cdots A_iW^{(i)}_{\gamma^p}.
\end{equation}
Therefore
\begin{equation}\label{eq:word-est}
 \sum_{\gamma\in\N_0^{m-1}}\norm{W_{(p,\gamma)}}
 \le\norm{A_i}^p L_i^{p+1}.
\end{equation}
Define the finite constant
\begin{equation}\label{eq:Ki}
 K_i:=\sum_{p=0}^{N-1}\norm{A_i}^pL_i^{p+1}<\infty.
\end{equation}
For any $\alpha$ with $\alpha_i<N$, and any
$\beta\in R_k(\alpha)$, one has $0\le\beta_i\le\alpha_i-1<N$.
Thus \eqref{eq:word-est} implies
\[
 \sum_{\beta\in R_k(\alpha)}\norm{W_\beta}\le K_i
 \qquad(1\le k\le m).
\]
The solution formula \eqref{eq:representation} yields
\begin{equation}\label{eq:slab}
 \norm{x_\alpha}\le
 K_i\left(\sum_{k=1}^m\norm{A_k}\right)\norm{\xi}_\infty
 \quad\text{if }\alpha_i<N.
\end{equation}
Combining these $m$ slab bounds with the uniform tail estimate
\eqref{eq:tail}, we obtain a constant $C>0$ for which
\begin{equation}\label{eq:global}
 \norm{x_\alpha}\le C\norm{\xi}_\infty
 \quad\text{for all }\alpha\in\N^m.
\end{equation}
This is precisely the required Lyapunov-stability bound.
\end{proof}

\begin{cor}\label{cor:equiv}
For \eqref{eq:FM} in the bounded-boundary/Pringsheim framework,
attractivity and asymptotic stability are equivalent.
\end{cor}
\begin{proof}
Attractivity implies Lyapunov stability by
Theorem~\ref{thm:main}; the reverse implication follows from the
definition of asymptotic stability.
\end{proof}

An immediate application combines our new implication with the
independently established spectral theorem.
\begin{cor}[Spectral test for attractivity]\label{cor:spectral}
The system \eqref{eq:FM} is attractive if and only if
\[
 \max_{(z_1,\ldots,z_m)\in\T^m}
 \rho\!\left(\sum_{k=1}^m z_kA_k\right)<1.
\]
\end{cor}

\begin{proof}
Combine Corollary~\ref{cor:equiv} and
Theorem~\ref{thm:FMcriterion}.
\end{proof}

\begin{cor}[Exponential interior decay]\label{cor:decay}
If \eqref{eq:FM} is attractive, then some $C>0$ and
$q\in(0,1)$ satisfy
\begin{equation}\label{eq:decay}
 \norm{x_\alpha}\le Cq^{\min_k\alpha_k}\norm{\xi}_\infty
 \quad\text{for all }\alpha\in\N^m.
\end{equation}
\end{cor}

\begin{proof}
By Corollary~\ref{cor:equiv}, the system is asymptotically stable.
Lemma~\ref{lem:summability} yields
$\norm{W_\beta}\le Mq^{\abs\beta}$ for some $q\in(0,1)$.
For $\beta\in R_k(\alpha)$, $\beta_k=\alpha_k-1$, so
\[
 \sum_{\beta\in R_k(\alpha)}\norm{W_\beta}
 \le Mq^{\alpha_k-1}(1-q)^{-(m-1)}.
\]
Insert this in \eqref{eq:representation} and absorb
$Mq^{-1}(1-q)^{-(m-1)}\sum_k\norm{A_k}$ into $C$.
\end{proof}

\section{Comparison with radial attractivity}\label{sec:comparison}
Bachelier et al. \cite{BachelierECC} use $c_0$ axis data and radial
convergence ($i+j\to\infty$) in dimension two, rather than bounded
face data and Pringsheim convergence.

\begin{ex}[Scalar averaging]\label{ex:average}
For
\begin{equation}\label{eq:average}
 x_{i,j}=\tfrac12(x_{i-1,j}+x_{i,j-1}),\qquad i,j\ge1,
\end{equation}
$\abs{x_{i,j}}\le\norm{\xi}_\infty$, so the system is Lyapunov
stable. Constant axis data equal to $1$ give $x_{i,j}=1$; thus it is
not attractive for arbitrary bounded data.
 
In contrast, with $c_0$ axis data, \eqref{eq:average} \emph{is} attractive in the radial sense $i+j\rightarrow\infty$.
For fixed
$k\ge1$, the coefficient of $x_{k,0}$ in $x_{i,j}$ is
\begin{equation}\label{eq:coefficient}
 2^{-(i-k+j)}\binom{i-k+j-1}{i-k}\quad(i\ge k,\ j\ge1),
\end{equation}
zero for $i<k$, and similarly for $x_{0,k}$. The bound
$2^{-p}\max_{0\le h\le p}\binom ph=O((p+1)^{-1/2})$ shows that each such
coefficient vanishes as $i+j\to\infty$. Hence finite-support data
give radial decay; uniform approximation and the stability bound extend
this to all $c_0$ data.
\end{ex}

Thus the same recurrence is radially attractive for $c_0$ data but
Lyapunov stable and non-attractive for arbitrary bounded data.

\section{Conclusion}
We proved that Pringsheim attractivity of a multidimensional FM system,
with arbitrary bounded data on the open coordinate faces, implies
Lyapunov stability. The proof supplements the Baire-category tail bound
with uniform control of the coordinate slabs, obtained by induction and
summability of word sums for lower-dimensional subsystems. When combined
with the spectral theorem of \cite{BourgetMcMillenPreprint}, this also
yields a spectral characterization of attractivity and a uniform
exponential decay bound in the minimum coordinate. Establishing a
parallel statement for multidimensional Roesser-type systems under the
same boundary and convergence conventions would be a natural next step.

\end{document}